\documentclass[a4,12pt]{article}
\usepackage{amsmath}
\usepackage{amssymb}

\begin{document}
\noindent
\begin{center}
  {\LARGE \bf A Rigidity Theorem for Affine
K$\ddot{a}$hler-Ricci Flat Graph }
  \end{center}
\noindent
\begin{center}
\begin{center}
   {\large An-Min Li\ \ and\ \  Ruiwei Xu  }\footnote{The first
   author is
   partially supported by NKBRPC(2006CB805905), NSFC 10631050 and RFDP.}\\[5pt]
  \end{center}
\end{center}

\noindent {\bf Abstract:} It is shown that any smooth strictly
convex global solution on $\mathbb{R}^n$ of
$$\det\left(\frac{\partial^{2}u}{\partial \xi_{i}\partial
\xi_{j}}\right) = \exp \left\{-\sum_{i=1}^n d_i \frac{\partial
u}{\partial \xi_{i}} - d_0\right\},$$ where $d_0$, $d_1$,...,$d_n$
are constants, must be a quadratic polynomial. This extends a
well-known theorem of J\"{o}rgens-Calabi-Pogorelov.

\vskip 6pt\noindent {\bf 2000 AMS Classification:} 53A15. \vskip
0.1pt\noindent {\bf Key words:} Pogorelov Theorem;
K$\ddot{a}$hler-Ricci Flat; Monge-Amp\`{e}re equation. \vskip 0.1pt

\section*{\S1. Introduction}
\vskip 0.1in\noindent

A well-known theorem of J\"{o}rgens $(n=2$ [J]), Calabi $(n\leq 5$
[Ca]), and Pogorelov $(n\geq 2$ [P]) states that any smooth strictly
convex solution of
$$\det\left(\frac{\partial^2f}{\partial x_i\partial x_j}\right)
=1\quad\quad\hbox{on}\ \ \mathbb{R}^n \leqno(1.1)$$ must be a
quadratic polynomial. In [C-Y] Cheng and Yau gave an analytical
proof. Recently Caffarelli and Li [C-L] extended the result for
classical solution to viscosity solution.

In this paper we study the following PDE
 $$\frac{\partial^{2}}{\partial x_{i}\partial x_{j}}\left(\log \det
\left( \frac{\partial^{2}f}{\partial x_{k}\partial x_{l}}
\right)\right)  = 0,\leqno(1.2)$$ or
$$\det\left(\frac{\partial^{2}f}{\partial x_{i}\partial
x_{j}}\right) = \exp \left\{\sum_{i=1}^n d_i x_i +
d_0\right\},\leqno(1.3)$$ where $d_0$, $d_1$,...,$d_n$ are
constants. Obviously, all solutions of $(1.1)$ satisfy (1.2).
Introduce the Legendre transformation of $f$
$$  \xi_{i}=\frac{\partial f}{\partial x_{i}},\quad i=
1,2,\dots,n,$$$$u(\xi_{1},\dots,\xi_{n})=  \sum_{i=1}^n
 x_{i}\frac{\partial f}{\partial x_{i}}- f(x). $$ In
terms of $\xi_1,...,\xi_n, u(\xi_1,...,\xi_n)$, the PDE (1.3) can
be written as
$$\det\left(\frac{\partial^{2}u}{\partial \xi_{i}\partial
\xi_{j}}\right) = \exp \left\{-\sum_{i=1}^n d_i \frac{\partial
u}{\partial \xi_{i}} - d_0\right\}.\leqno(1.4)$$ Note that, under
the Legendre transformation, the PDE (1.1) reads
$$\det\left(\frac{\partial^2u}{\partial \xi_i\partial \xi_j}\right)
=1\quad\quad\hbox{on}\ \ \mathbb{R}^n .\leqno(1.1)'$$  Given any
smooth, bounded convex domain $\Omega \subset \mathbb{R}^n$ and
any smooth boundary value $\phi$, the existence of the solution of
the boundary problem
$$\det\left(\frac{\partial^{2}u}{\partial \xi_{i}\partial
\xi_{j}}\right) = \exp \left\{-\sum_{i=1}^n d_i \frac{\partial
u}{\partial \xi_{i}} -
d_0\right\}\;\;\;\hbox{in}\;\;\Omega,\;\;\;u = \phi\;\;\;\hbox{on}
\;\;\;\partial \Omega\leqno(1.4)'$$ is well-known. So there are
many locally solutions to the PDE $(1.4)$. In this paper we prove
the following theorem \vskip0.1in\noindent {\bf Main Theorem.}
{\it Let $u(\xi_1,...,\xi_n)$ be a $C^{\infty}$ strictly convex
function defined on whole $\mathbb{R}^n$. If $u(\xi)$ satisfies
the PDE (1.4) , then $u$ must be a quadratic polynomial.} \vskip
0.1in\noindent

\vskip 0.1in The PDE (1.2) arises naturally in the construction of
Ricci flat K$\ddot{a}$hler-affine metric for affine manifolds.  An
affine manifold is a manifold which can be covered by coordinate
charts so that the coordinate transformations are given by
invertible affine transformations. Let $M$ be an affine manifold.
A K$\ddot{a}$hler affine metric or Hessian metric $G$ on $M$ is a
Riemannian metric on $M$ such that locally, for affine coordinates
$(x_{1},x_{2},\cdots, x_{n})$, there is a potential $f$ such
that\vskip 0.1in \noindent
$$G_{ij}= \frac{\partial^2f}{\partial x_i\partial x_j} .$$
The pair $\left(M, G\right)$ is called a K$\ddot{a}$hler affine
manifold or a Hessian manifold, and $G$ is called K$\ddot{a}$hler
affine metric. K$\ddot{a}$hler affine metric was first studied by
Cheng and Yau in [C-Y-1]. For more details about Hessian manifolds
please see [Sh]. Following Cheng and Yau we introduce the concepts
of the K$\ddot{a}$hler Ricci curvature and the K$\ddot{a}$hler
scalar curvature of $G$ on $M$. It is easy to see that the tangent
bundle $TM$ is a complex manifold with a natural complex structure
in the following way. For coordinate chart $(x_1,x_2,\dots,x_n)$,
we can consider a tube over the coordinate neighborhood with
complex coordinate system $(x_1+iy_1,x_2+iy_2,\dots,x_n+iy_n)$.
The Hessian metric $G$ was naturally extended to be a
K$\ddot{a}$hler metric of the complex manifold $TM$. The Ricci
curvature tensor and the scalar curvature of this K$\ddot{a}$hler
metric are given by respectively
$$ R_{ij}= - \frac{\partial^{2}}{\partial x_{i}\partial x_{j}}
\left(\log \det\left(f_{kl}\right)\right),\quad
R=-\frac{1}{2}\sum^{n}_{i=1}\sum^{n}_{j=1}f^{ij}\frac{\partial^{2}(\log
\det(f_{kl}))} {\partial x_{i}\partial x_{j}}. $$ It is obvious
that the restrictions of $R_{ij}$ and $R$ to $M$ are tensors of
$M$. We also call $R_{ij}$ and $R$ the K$\ddot{a}$hler Ricci
curvature and the K$\ddot{a}$hler scalar curvature of $G$ on $M$.
We say that the metric $G$ is K$\ddot{a}$hler-Ricci flat if (1.2)
holds on $M$ everywhere. In this geometric language, our Main
Theorem can be stated as

\vskip 0.1in\noindent{\bf Main Theorem.} {\it Let $M$ be a graph
given by a smooth strictly convex function $x_{n+1} =
f(x_1,...,x_n)$ defined in a domain $\Omega$. If the Hessian
metric of $M$ is K$\ddot{a}$hler-Ricci flat and the image of $M$
under the normal mapping is whole $\mathbb{R}^n$, then $f$ must be
a quadric.} \vskip 0.1in\noindent {\bf Remark 1.} In [J-L] the
authors have proved that

\vskip 0.1in\noindent {\bf Theorem.} {\it Let $M$ be a
K$\ddot{a}$hler affine manifold. If the Hessian metric of $M$ is
K$\ddot{a}$hler-Ricci flat and complete, then $M$ must be
$\mathbb{R}^n/\Gamma$, where $\Gamma$ is a subgroup of isometries
which acts freely and properly discontinuously on $\mathbb{R}^n$.}

 \vskip 0.1in\noindent
{\bf Remark 2.} From our proof of the Main Theorem the following
stronger version is also true: \vskip 0.1in\noindent {\bf Main
Theorem'.} {\it Let $u(\xi_1,...,\xi_n)$ be a $C^{\infty}$
strictly convex function defined in a convex domain $\Omega\subset
\mathbb{R}^n.$ If $u(\xi)$ satisfies the PDE (1.4) and if
$u(p)\rightarrow \infty$ as $p\rightarrow \partial \Omega$, then
$u$ must be a quadratic polynomial.} \vskip 0.1in\noindent {\bf
Remark 3.} The global solution of the PDE (1.3) on the
$x-coordinate$ plane $\mathbb{R}^n$ is not unique. For example,
$$ f(x_1,...,x_n)=\sum _{i=1}^n x_i^2,\quad \hbox{and}\quad
f(x_1,...,x_n) = \exp\{x_1\} + \sum _{i=2}^n x_i^2$$ are global
solutions of the PDE (1.3).

 \vskip 0.1in\noindent
{\bf Remark 4.} Our study in this paper is based on the following
differential inequality for $\Phi$ (for details see Proposition
3.1 below)
$$\Delta \Phi\geq\frac{n}{n-1}\frac{\left\|\nabla \Phi\right\|^2}{\Phi}+
\frac{n^{2}-3n-10}{2(n-1)}\langle\nabla \Phi,\nabla \log
\rho\rangle + \frac{(n+2)^{2}}{n-1}\Phi^{2}.$$  This type of
differential inequality for $\Phi$ first appeared in [L-J-1], in
which Li and Jia  announced that they solved the Chern's
conjecture for 2-dimension and 3-dimension. While Trudinger and
Wang solved Chern's conjecture for 2-dimension in [T-W]. Li and
Jia's method, which is quite different from that of Trudinger and
Wang, is to estimate $\Phi$ and $\|\nabla f\|$ based on the
differential inequality:
$$\Delta^B \Phi \geq
\frac{n}{2(n-1)}\frac{\left\|\nabla \Phi\right\|^2_{G^B}}{\Phi} -
\frac{n^2-n-2}{2(n-1)}\langle\nabla \Phi,\nabla \log
\rho\rangle_{G^B}$$
$$ + \left( 2 - \frac{(n-2)^2(n-1)}{8n} - \frac{n^2
-2}{2(n-1)}\right)
 \frac{\Phi^2}{\rho},$$
where $G^B$ is the Blaschke metric and $\Delta^B$ is the Laplacion
with respect to $G^B$.

\vskip 0.1in However, Li later found  a gap in their proof, so the
full research paper is not published. In [L-J-2] the author use
the similar differential inequality to prove Bernstein properties
for some more general fourth order nonlinear PDE for 2 dimension.
As a corollary, they fix the gap to 2 dimensional Chern'
conjecture. So far the 3 dimensional Chern' conjecture is open.

\section*{\S2. Preliminaries }
\vskip 0.1in\noindent

Let $f(x_1,...,x_n)$ be a $C^{\infty}$ strictly convex function
defined on a domain $\Omega \subset \mathbb{R}^n$. Denote
$$M:=\{(x,f(x))|x_{n+1}=f(x_1,...,x_n),\;\;\;(x_1,...,x_n)
\in \Omega\}.$$  We choose the canonical relative normalization $Y
=(0,0,...,1)$. Then, in terms of language of the relative affine
differential geometry, $G$ is the relative metric with respect to
the normalization $Y$. Denote by $y = (x_1,...,x_n,
f(x_1,...,x_n))$, the position vector of $M$. We have
$$y_{,ij}= \sum A_{ij}^ky_k + f_{ij}Y.\leqno(2.1)$$
The conormal field $U$ is given by
$$U = \left( - f_1,...,-f_n, 1\right).\leqno(2.2)$$
We recall some fundamental formulas for the graph $M$ without proof,
for details see [P-1]. The Levi-Civita connection with respect to
the metric $G$ is
$$\Gamma^k_{ij} = \frac{1}{2}\sum f^{kl}f_{ijl},\leqno(2.3)$$
The Fubini-Pick tensor $A_{ijk}$ and the Weingarten tensor are given
by
$$A_{ijk} = - \frac{1}{2}f_{ijk},\;\;\;\;B_{ij} = 0.\leqno (2.4)$$
The relative Pick invariant is
$$J = \frac{1}{4n(n-1)}\sum f^{il}f^{jm}f^{kn}f_{ijk}f_{lmn}.\leqno(2.5)$$
The Gauss equations and the Codazzi equations read
$$R_{ijkl}=\sum f^{m h}( A_{jkm}A_{hil}-A_{ikm}A_{hjl}), \leqno(2.6)$$
$$A_{ijk,l} = A_{ijl,k}.\leqno(2.7)$$
From (2.6) we have
$$R_{ik}= \sum f^{mh}f^{lj}(A_{iml}A_{hjk} - A_{imk}A_{hlj}).\leqno (2.8)$$
Denote
$$\rho = \left [\det (f_{ij})\right ]^{-\frac{1}{n+2}},\quad
\Phi = \frac{\left\|\nabla \rho \right\|^2}{\rho ^2}.\leqno
(2.9)$$ Let $\Delta $ be the laplacian with respect to the Calabi
metric, which is defined by
$$\Delta = \frac{1}{\sqrt{\det(G_{kl})}}\sum
\frac{\partial}{\partial x_i}
\left(G^{ij}\sqrt{\det(G_{kl})}\frac{\partial}{\partial
x_j}\right).\leqno(2.10)$$ By a direct calculation from (2.10) we
have
$$ \Delta =\sum f^{ij}\frac{\partial^2}{\partial x_i
\partial x_j} + \frac{n+2}{2}\frac{1}{\rho} \sum
f^{ij}\frac{\partial \rho}{\partial x_j} \frac{\partial }
{\partial x_i}\leqno(2.11) $$$$ =  \sum
u^{ij}\frac{\partial^2}{\partial \xi_i
\partial \xi_j} - \frac{n+2}{2}\frac{1}{\rho } \sum
u^{ij}\frac{\partial \rho}{\partial \xi_j} \frac{\partial }
{\partial \xi_i},  $$
$$ \Delta f   = n + \frac{n+2}{2}\frac{1}{\rho}
\langle\nabla \rho, \nabla f\rangle,\leqno(2.12)$$$$ \Delta u   =n
- \frac{n+2}{2}\frac{1}{\rho} \langle\nabla \rho, \nabla
u\rangle.\leqno(2.13) $$

\section*{\S3.  Calculation of $\Delta \Phi$ } 
\vskip 0.1in\noindent

The following proposition is proved in [J-L], however, we include
here for the reader's convenience. \vskip 0.1in\noindent

{\bf Proposition 3.1} {\it Let $f(x_1,...,x_n)$ be a $C^{\infty}$
strictly convex function satisfying the PDE (1.3). Then the
following estimate holds
$$\Delta \Phi\geq\frac{n}{n-1}\frac{\left\|\nabla \Phi\right\|^2}{\Phi}+
\frac{n^{2}-3n-10}{2(n-1)}\langle\nabla \Phi,\nabla \log \rho\rangle
+ \frac{(n+2)^{2}}{n-1}\Phi^{2}.$$} \vskip 0.1in\noindent

{\bf Proof.} From the PDE (1.4) we have
$$ 0 = \frac{\partial^{2}}{\partial x_{i}\partial x_{j}}
\left(\log \det\left(f_{kl}\right)\right)=
-(n+2)\left(\frac{\rho_{ij}}{\rho}-
\frac{\rho_{i}}{\rho}\frac{\rho_{j}}{\rho}\right),\leqno(3.1)$$
where $\rho_{i}=\frac{\partial\rho}{\partial x_{i}}$ and $\rho_{ij}=
\frac{\partial^{2}\rho}{\partial x_{i}\partial x_{j}}.$ It follows
that
$$\Delta \rho=\frac{n+4}{2}\frac{\left\|\nabla\rho\right\|^{2}}{\rho}.
\leqno(3.2)$$ Let $p\in M$, we choose a local orthonormal frame
field of the metric $G$ around $p$. Then
$$\Phi=\frac{\sum (\rho_{,j})^2}{\rho^{2}},\quad
\Phi_{,i}=2\sum\frac{\rho_{,j}\rho_{,ji}}{\rho^{2}}-2\rho_{,i}
\frac{\sum (\rho_{,j})^2}{\rho^{3}},$$
$$\Delta \Phi=2\frac{\sum (\rho_{,ji})^2}{\rho^{2}}+2\sum\frac{\rho_{,j}
\rho_{,jii}}{\rho^{2}}-8\sum\frac{\rho_{,j}\rho_{,i}\rho_{,ji}}
{\rho^{3}}-(n-2)\frac{\left(\sum
(\rho_{,j})^2\right)^{2}}{\rho^{4}},$$ where we used (3.2). In the
case $\Phi(p)=0$, it is easy to get, at $p,$ $$\Delta \Phi\geq
2\frac{\sum (\rho_{,ij})^2}{\rho^{2}}.$$ Now we assume that $\Phi
(p)\neq 0$. Choose a local orthonormal frame field of the metric
$G$ around $p$ such that
$\rho_{,1}(p)=\left\|\nabla\rho\right\|(p)>0,\;\;\rho_{,i}(p)=0$
for all $i>1$. Then
$$\Delta \Phi=2\frac{\sum (\rho_{,ij})^2}{\rho^{2}}+2\sum\frac{\rho_{,j}
\rho_{,jii}}{\rho^{2}}-8\frac{(\rho_{,1})^2\rho_{,11}}
{\rho^{3}}-(n-2)\frac{(\rho_{,1})^4}{\rho^{4}}. \leqno(3.3)$$
Applying an elementary inequality
$$a_1^2+a_2^2+\dots+a_{n-1}^2\geq\frac{(a_1+a_2+\dots+a_{n-1})^2}{n-1}$$
and (3.2), we obtain
$$2\frac{\sum (\rho_{,ij})^2}{\rho^{2}}\geq
2\frac{(\rho_{,11})^2}{\rho^{2}}+4
\frac{\sum_{i>1}(\rho_{,1i})^2}{\rho^{2}}+2
\frac{\sum_{i>1}(\rho_{,ii})^2}{\rho^{2}} \geq
2\frac{(\rho_{,11})^2}{\rho^{2}}+
4\frac{\sum_{i>1}(\rho_{,1i})^2}{\rho^{2}}\leqno(3.4)$$$$
+\frac{2}{n-1} \frac{(\Delta\rho -\rho_{,11})^2}{\rho^{2}}
\geq\frac{2n}{n-1}\frac{(\rho_{,11})^2}
{\rho^{2}}+4\frac{\sum_{i>1}(\rho_{,1i})^2}{\rho^{2}}-2\frac{n+4}{n-1}
\frac{(\rho_{,1})^2\rho_{,11}}{\rho^{3}}+\frac{(n+4)^{2}}{2(n-1)}
\frac{(\rho_{,1})^4}{\rho^{4}}.$$ An application of the Ricci
identity shows that\vskip 0.1in \noindent
$$\frac{2}{\rho^{2}}\sum\rho_{,j}\rho_{,jii}=\frac{2}{\rho^2}
(\Delta\rho)_{,1}\rho_{,1}+2R_{11}\frac{(\rho_{,1})^2}{\rho^2}
\leqno(3.5)$$$$=2(n+4)\frac{(\rho_{,1})^2\rho_{,11}}{\rho^{3}}-(n+4)\frac{(\rho_{,1})^4}
{\rho^{4}}+2R_{11}\frac{(\rho_{,1})^2}{\rho^{2}}.$$ Substituting
(3.4) and (3.5) into (3.3) we obtain
$$ \Delta \Phi\geq  \frac{2n}{n-1}\frac{(\rho_{,11})^2}{\rho^{2}}
+\left(2n-2\frac{n+4}{n-1}\right)\frac{(\rho_{,1})^2
\rho_{,11}}{\rho^{3}}+2R_{11}\frac{(\rho_{,1})^2}{\rho^{2}}
 \leqno (3.6)$$$$+\left(\frac{ (n+4)^{2}}{2(n-1)}-2(n+1)\right)
\frac{(\rho_{,1})^4}{\rho^{4}}+ 4
\frac{\sum_{i>1}(\rho_{,1i})^2}{\rho^2}. $$  Note that
$$\sum\frac{(\Phi_{,i})^2}{\Phi}=4\frac{\sum (\rho_{,1i})^2}
{\rho^{2}}-8\frac{(\rho_{,1})^2\rho_{,11}}{\rho^{3}}
+4\frac{(\rho_{,1})^4}{\rho^{4}}.\leqno(3.7)$$ Then (3.6) and
(3.7) together give us
$$\Delta \Phi\geq
\frac{n}{2(n-1)}\frac{\sum (\Phi_{,i})^2}{\Phi}+
\left(\frac{2n-8}{n-1}+2n \right)
\frac{(\rho_{,1})^2\rho_{,11}}{\rho^{3}}\leqno(3.8)$$
$$+2R_{11}\frac{(\rho_{,1})^2}{\rho^{2}}+\left[\frac{(n+4)^{2}}
{2(n-1)}-2(n+1)-\frac{2n}{n-1}\right]\frac{(\rho_{,1})^4}
{\rho^{4}}.$$ From (3.1) we easily obtain
$$\rho_{,ij}=\rho_{ij} + A_{ij1}\rho_{,1}=
\frac{\rho_{,i}\rho_{,j}}{\rho} + A_{ij1}\rho_{,1}.$$ Thus we get
$$ \Phi_{,i}=\frac{2\rho_{,1}\rho_{,1i}}{\rho^2}
-2\frac{\rho_{,i}(\rho_{,1})^2}{\rho^3}=
2A_{i11}\frac{(\rho_{,1})^2}{\rho^2},\quad\sum
\Phi_{,i}\frac{\rho_{,i}}{\rho}=2\frac{(\rho_{,1})^2\rho_{,11}}{\rho^3}
-2\frac{(\rho_{,1})^4}{\rho^4},\leqno(3.9)$$$$ \quad \frac{\sum
(\Phi_{,i})^2}{\Phi} = 4\sum
(A_{i11})^2\frac{(\rho_{,1})^2}{\rho^2},\;\;\;\; \sum
\Phi_{,i}\frac{\rho_{,i}}{\rho} =
2A_{111}\frac{(\rho_{,1})^3}{\rho^3}.\leqno(3.10)$$ By the same
method as deriving (3.4) we get
$$ \sum (A_{ml1})^2\geq (A_{111})^2+2\sum_{i>1}(A_{i11})^2+\sum_{i>1} (A_{ii1})^2
\leqno(3.11)$$$$\geq (A_{111})^2+2\sum_{i>1}(A_{i11})^2
+\frac{1}{n-1}\left(\sum A_{ii1}-A_{111}\right)^2$$$$
\geq\frac{n}{n-1}\sum (A_{i11})^2- \frac{2}{n-1}A_{111} \sum
A_{ii1} + \frac{1}{n-1}\left(\sum A_{ii1}\right)^2.$$ Therefore,
by (2.8), (3.10) and (3.11), we obtain
$$2R_{11}\frac{(\rho_{,1})^2}{\rho^{2}}=2\sum (A_{kj1})^2\frac{(\rho_{,1})^2}{\rho^2}
-(n+2)A_{111}\frac{(\rho_{,1})^3}{\rho^3}\leqno(3.12)$$$$\geq\frac{n}{2(n-1)}
\frac{\sum (\Phi_{,i})^2}{\Phi}-\frac{(n+2)(n+1)}{2(n-1)}
\sum\Phi_{,i}\frac{\rho_{,i}}{\rho}+\frac{(n+2)^{2}}{2(n-1)}
\frac{(\rho_{,1})^4}{\rho^{4}}. $$ Then inserting (3.12) and (3.9)
into (3.8) we have
$$\Delta \Phi\geq\frac{n}{n-1}\frac{\sum (\Phi_{,i})^2}{\Phi}+
\frac{n^{2}-3n-10}{2(n-1)}\sum \Phi_{,i}\frac{\rho_{,i}}{\rho}+
\frac{(n+2)^{2}}{n-1}\Phi^{2}.\quad\Box\leqno(3.13)$$ \vskip 0.1in

\section*{\S4. Proof of Main Theorem for $n\leq 4$}
\vskip 0.1in\noindent

In the case $n\leq 4$ the proof of the Main Theorem is relative
simple, we first consider this case.

We shall show that  $\Phi = 0$ on $M$ everywhere,  namely,
$\det\left(\frac{\partial^{2}u}{\partial
\xi_{i}\partial\xi_{j}}\right)=const.$ Therefore the main Theorem
follows by J-C-P Theorem. By a coordinate translation transformation
and by subtracting a linear function we may suppose that
$$u(0)= 0,\;\;\; u(\xi)\geq 0.$$
Then for any $C>0$ the set
$$\bar{S}_u(0,C) := \{ \xi \in \mathbb{R}^n | u(\xi) \leq C\}$$
is compact. Consider the function
$$L = \exp \left \{-\frac{m}{C-u}\right \}\Phi$$
defined on $\bar{S}_u(0,C)$, where $m$ is a positive constant to be
determined later. Clearly, $L$ attains its supremum at some interior
point $p^*$. Then, at $p^*$,\vskip 0.1in \noindent
$$ \frac{\Phi_{,i}}{\Phi} -  hu_{,i}= 0,\leqno(4.1)$$$$
\frac{\Delta \Phi}{\Phi}- \frac{\sum (\Phi_{,i})^2}{\Phi^2} -
h'\sum (u_{,i})^2 - h \Delta u \leq 0,\leqno(4.2) $$ where and
later we denote
$$h=\frac{m}{(C-u)^2},\;\;\;h'=\frac{2m}{(C-u)^3},$$
and "," denotes the covariant derivatives with respect to the metric
$G$. Inserting (3.13) (2.13) and (4.1) into (4.2) we get
$$\frac{(n+2)^{2}}{n-1}\Phi + \left(\frac{1}{n-1}h^2 - h'\right)\sum (u_{,i})^2
- nh  + \frac{(n+2)(n-3)}{(n-1)}h \frac{\sum
\rho_{,i}u_{,i}}{\rho}\leq 0.\leqno(4.3)$$ By the Schwarz's
inequality
$$\frac{(n+2)(n-3)}{(n-1)}h \frac{\sum \rho_{,i}u_{,i}}{\rho}\leq
\frac{1}{2(n-1)}h^2 \sum
(u_{,i})^2+\frac{(n+2)^2(n-3)^2}{2(n-1)}\Phi.$$ Therefore
$$\frac{(n+2)^{2}(2-(n-3)^2)}{2(n-1)}\Phi + \left(\frac{1}{2(n-1)}h^2
 - h'\right)\sum (u_{,i})^2- nh \leq 0.\leqno(4.4)$$
In the case $n\leq 4$ we have
$$\frac{(n+2)^{2}}{2(n-1)}\Phi + \left(\frac{1}{2(n-1)}h^2 -
h'\right)\sum (u_{,i})^2- nh \leq 0.\leqno(4.5)$$ We choose $m =
8(n-1)C$, then $\frac{1}{2(n-1)}h^2 - h' \geq 0$. It follows that,
at $p^*$,
$$\exp  \left\{-\frac{8(n-1)C}{C-u}\right\}\Phi
\leq n \exp\left\{-\frac{m}{C-u}\right\}h\leq\frac{b}{C},
\leqno(4.6)
$$ where $b$ is a constant depending only on $n$. In the
calculation of (4.6) and later we often use the fact that
$\exp\left\{-\frac{m}{C-u}\right\}\frac{m^2}{(C-u)^2} $ has a
universal upper bound. Since $L$ attains its supremum at $p^*$,
(4.6) holds everywhere in $\bar{S}_u(0,C)$. For any fixed point
$p$, we let $C\rightarrow \infty$ then $\Phi(p) = 0$. Therefore
$\Phi = 0$ everywhere on $M$. \quad $\Box$

\section*{\S5. Estimate for $\sum \left(\frac{\partial u}{\partial
\xi_i}\right)^2$} \vskip 0.1in \noindent

For general dimensions $(n>4)$ the proof of the Main Theorem is much
more difficult than $n\leq 4$, it needs more estimates. In this
section we estimate $\sum \left(\frac{\partial u}{\partial
\xi_i}\right)^2$. Let $\Omega \subset \mathbb{R}^n$ be a bounded
convex domain. It is well-known (see [G]) that there exists a unique
ellipsoid $E$, which attains the minimum volume among all the
ellipsoids that contain $\Omega $ and that are centered at the
center of mass of $\Omega $, such that
$$n^{-\frac{3}{2}} E \subset \Omega  \subset E,$$ where
$n^{-\frac{3}{2}} E$ means the $n^{-\frac{3}{2}}$ -dilation of $E$
with respect to its center. Let $T$ be an affine transformation such
that $T(E)=B(0,1)$, the unit ball. Put $\tilde{\Omega}= T(\Omega)$.
Then
$$B(0,n^{-\frac{3}{2}}) \subset \tilde{\Omega} \subset
B(0,1).\leqno(5.1)$$ A convex domain $\Omega $ is called normalized
if it satisfies (5.1). Let $u$ be a smooth strictly convex function
defined on $\Omega$ such that
$$\inf_{\Omega}u(\xi)= u(p) = 0,\;\;u|_{\partial \Omega}
 =1.\leqno(5.2)$$
A strictly convex function defined on $\Omega $ is called normalized
at $p$ if (5.2) holds. \vskip 0.1in \noindent {\bf Lemma 5.1} {\it
Let $\Omega_k$ be a sequence of smooth and normalized convex
domains, $u^{(k)}$ be a sequence of smooth strictly convex functions
defined on $\Omega_k$, normalized at $p_{k}$. Then there are
constants $d>1$, $b>0$ independent of $k$ such that
$$\frac{\sum_i(\frac {\partial u^{(k)}}{\partial \xi_i})^2 }
{(d+f^{(k)})^2}\leq b, \ \ \ \ k=1,2,\dots \quad\hbox{on}\ \
\bar{\Omega}_k.$$ }

{\bf Proof.} We may suppose by taking subsequence that $\Omega_k$
converges to a convex domain $\Omega $ and $u^{(k)}$ converges to a
convex function $u^\infty $, locally uniformly in $\Omega$.
Obviously, we have the uniform estimate
$$\sum \left(\frac{\partial u^{(k)}}{\partial \xi_i}\right)^2(0)\leq
4n^3.\leqno(5.3)$$ For any $k$, let
$$\tilde{u}^{(k)} = u^{(k)} - \sum \frac{\partial
u^{(k)}}{\partial \xi_i}(0)\xi_i - u^{(k)}(0).\leqno(5.4)$$ Then
$$\tilde{u}^{(k)}(0) = 0,\;\;\;\;\tilde{u}^{(k)}(\xi)\geq 0,\;\;\;
\tilde{u}^{(k)}|_{\partial \Omega_k}\leq C_0,$$ where $C_0$ is a
constant depending only on $n$. As $B(0,n^{-\frac{3}{2}})
\subset\Omega_k$, we have
$$\frac{\mid\nabla \tilde{u}^{(k)}\mid^2}{(1+\tilde{f}^{(k)})^2}\leq
\mid\nabla \tilde{u}^{(k)}\mid^2\leq
\frac{C_0^2}{dist(B(0,2^{-1}n^{-\frac{3}{2}}),\partial\Omega_k)^2}\leq
4n^{3}C^2_0 \ \ \ \ \ \leqno(5.5)$$ on
$B(0,2^{-1}n^{-\frac{3}{2}})$, where $\tilde{f}^{(k)}$ is the
Legendre transformation of $\tilde{u}^{(k)}$ relative to $0$. For
any $p\in\bar{\Omega}_k\backslash B(0,2^{-1}n^{-\frac{3}{2}})$, we
may suppose that $p =(\xi_1,0,\dots,0)$ with $\xi_1> 0$ by an
orthonormal transformation. Then, at $p$,
$$C_0+\tilde{f}^{(k)}\geq\tilde{u}^{(k)}+\tilde{f}^{(k)}=
\frac{\partial \tilde{u}^{(k)}}{\partial\xi_1}\xi_1.$$ It follows
that
$$\frac{\left(\frac {\partial \tilde{u}^{(k)}}{\partial \xi_1}\right)^2 }
{(C_0+\tilde{f}^{(k)})^2}<\frac{1}{\xi_1^2}<4n^3.$$ Therefore
there exist constants $\tilde{d}>1,$ $\tilde{b}>0$ depending only
on $n$ such that
$$\frac{\left(\frac {\partial \tilde{u}^{(k)}}{\partial r}\right)^2 }
{(\tilde{d}+\tilde{f}^{(k)})^2}< \tilde{b},\leqno(5.6)$$ where
$\frac{\partial}{\partial r}$ denotes the radial derivative.
 Note that
$$ \frac{\partial\tilde{u}^{(k)}}{\partial \xi_i}= \frac{\partial
u^{(k)}}{\partial \xi_i}- \frac{\partial{u}^{(k)}}{\partial
\xi_i}(0),\ \ \ \tilde{f}^{(k)}= f^{(k)} + u^{(k)}(0).
 \leqno(5.7)$$
It follows from (5.3) and (5.4) that
$$\left(\frac {\partial u^{(k)}}{\partial r}\right)^2 \leq
2\left(\frac {\partial \tilde{u}^{(k)}}{\partial r}\right)^2 +
8n^3.$$ Then
$$\frac{\left(\frac {\partial u^{(k)}}{\partial r}\right)^2 }
{(d'+ f^{(k)})^2}< b',\leqno(5.8)$$ for some constants $d'>1$,
$b'>0$ independent of $k$. Note that
$$|\nabla u^{(k)}(p)|=\frac{1}{\cos\alpha_k}\left|\frac
 {\partial u^{(k)}}{\partial r}(p)\right|,\leqno(5.9)$$
where $\alpha_k$ is the angle between vectors $\nabla u^{(k)}(p)$
and $\frac {\partial u^{(k)}}{\partial r }(p).$ Since $u^{(k)}= 1
$ on $\partial \Omega_k$, $\nabla u^{(k)}(p)$ is perpendicular to
the boundary of $\Omega_k$ at any $p\in \partial \Omega_k$. As
$\Omega$ is convex and $0\in \Omega$, it follows that
$\frac{1}{\cos\alpha_k}$ have a uniform upper bound. Then the
Lemma 5.1 follows. $\Box$
 \vskip 0.1in \noindent
{\bf Remark 5.2} We may choose $d$ in Lemma 5.1 such that the
following holds for any $k$
$$\frac{|u+f^{(k)}|}{d+f^{(k)}}\leq 1.\leqno(5.10)$$

\section*{\S6.  Estimates of $\rho$, $\rho^\alpha\Phi$ and $\sum u_{ii}$}
\vskip 0.1in \noindent

From now on we assume that $n\geq 5$. In this section we prove
some estimates which we need in the next section. Suppose that
$p\in \Omega$ and $u$ is normalized at $p$. For any positive
number $C\leq 1$, denote
$$S_u(p,C) = \left\{ \xi\in\Omega | u(\xi)<C \right \},\quad
\bar{S}_u(p,C) = \left\{ \xi\in\Omega | u(\xi)\leq C \right \}.$$
Introduce notations:
$$\mathcal {A}:=\max_{p\in S_u(p,C)}\left\{\exp\left\{-\frac{m}
{C-u}\right\}\frac{\rho^\alpha\Phi}{(d+f)^{\frac{2n\alpha}{n+2}}}
\right\},$$
$$\mathcal {B}:=\max_{p\in S_u(p,C)}\left\{\exp\left\{-\frac{m}
{C-u}+H\right\}\frac{(h+2\alpha)\rho^\alpha
}{(d+f)^{\frac{2n\alpha}{n+2}}}\right\},$$ where
$$\alpha =\frac{(n+2)(n-3)}{2} + \frac{n-1}{4} ,\quad m=32(n+2)C,
 \quad H=\epsilon\frac{\sum x_k^2}{(d+f)^2}.$$
From Lemma 5.1, we always choose small enough constant $\epsilon$
such that $H<\frac{1}{30}$ in this section.

We prove the following lemmas, which play important role in the
proof of the Main Theorem. \vskip 0.1in \noindent

{\bf Lemma 6.1} {\it  Let $u$ be a smooth and strictly convex
function defined in $\Omega$ which satisfies the equation (1.4).
Suppose that $u$ is normalized at $0$ and the section
$\bar{S}_u(p,C)$ is compact. And assume that  there are constants
$b_1\geq 0$, $d>1$ such that
$$\frac{\sum x_k^2}{(d+f)^2}\leq b_1$$ on $\bar{S}_u(p,C)$.
Then there is a constant $d_1>0$, depending only on $n$, $b_1$ and
$C$, such that $$\mathcal {A}\leq d_1,\quad \mathcal {B}\leq d_1.$$}

{\bf Proof.} Firstly, we show $\mathcal{A}\leq 10\mathcal {B}.$ To
this end, consider the following function
$$F = \exp \left\{-\frac{m}{C-u } \right\}\frac{ \rho^{\alpha}\Phi
}{(d + f)^{\frac{2n\alpha}{n+2}}} $$defined on $S_u(p,C)$.
Clearly, $F$ attains its supremum at some interior point $p^*$ of
$S_u(p,C)$. Thus, at $p^*$,
$$\frac{\Phi_{,i}}{\Phi} + \alpha \frac{\rho_{,i}}{\rho} -
\frac{2n\alpha}{n+2}\frac{f_{,i}}{d+f} - hu_{,i}= 0,\leqno(6.1)$$
$$\frac{\Delta \Phi}{\Phi}-
\frac{\sum (\Phi_{,i})^2}{\Phi^2} + \frac{n+2}{2}\alpha \Phi -
\frac{2n\alpha}{n+2}\frac{\Delta f}{d+f} +
\frac{2n\alpha}{n+2}\frac{\sum (f_{,i})^2}{(d+f)^2} - h'\sum
(u_{,i})^2 - h \Delta u \leq 0,\leqno(6.2)$$ where "," denotes the
covariant derivatives with respect to the metric $G$. In the
calculation of (6.2) we used (3.2). Inserting (2.12), (2.13) and
(3.13) into (6.2) we get
$$\left[\frac{(n+2)}{2}\alpha + \frac{(n+2)^{2}}{n-1}\right]\Phi +
\frac{1}{(n-1)}\frac{\sum (\Phi_{,i})^2}{\Phi^2} +
\frac{n+2}{2}h\frac{\sum u_{,i}\rho_{,i}}{\rho} - n\alpha
\frac{\sum f_{,i}\rho_{,i}}{(d+f)\rho}$$
$$+ \frac{2n\alpha}{n+2}\frac{\sum (f_{,i})^2}{(d+f)^2}+ \frac{n^2 - 3n -
10}{2(n-1)}\sum \frac{\Phi_{,i}}{\Phi}\frac{\rho_{,i}}{\rho} -
h'\sum (u_{,i})^2  - nh - \frac{2n\alpha}{n+2}\frac{n}{d+f}\leq
0.$$ Using (6.1) yields
$$\frac{1}{(n-1)}\sum \left[hu_{,i} +
\frac{2n\alpha}{n+2}\frac{f_{,i}}{d+f}- \alpha
\frac{\rho_{,i}}{\rho}\right]^2 + \left[\frac{2(n+2)}{n-1}\alpha +
\frac{(n+2)^{2}}{n-1}\right]\Phi \leqno(6.3)$$
$$+ \frac{(n+2)(n-3)}{n-1}h\frac{\sum u_{,i}\rho_{,i}}{\rho}
- \frac{4n\alpha}{n-1}\frac{\sum f_{,i}\rho_{,i}}{(d+f)\rho} +
\frac{2n\alpha}{n+2}\frac{\sum (f_{,i})^2}{(d+f)^2} - h'\sum
(u_{,i})^2 - nh - \frac{2n^2\alpha}{n+2}\leq 0.$$ Note that
$$\frac{|\sum u_{,i}f_{,i}|}{d+f} = \frac{|\sum \frac{\partial u}{\partial
\xi_i}\frac{\partial f}{\partial x_k}u_{kj}u^{ij}|}{d+f} =
\frac{|\sum \xi_i\frac{\partial u}{\partial \xi_i}|}{d+f} =\frac{
|u + f|}{d+f}\leq 1.\leqno(6.4)$$ Inserting (6.4) into (6.3), we
have
$$\frac{1}{(n-1)}h^2\sum (u_{,i})^2 + \left[
\frac{4n^2\alpha^2}{(n+2)^2(n-1)}+ \frac{2n\alpha}{n+2}\right]
\frac{\sum (f_{,i})^2}{(d+f)^2}- \frac{1}{2} h\frac{\sum
u_{,i}\rho_{,i}}{\rho}\leqno(6.5)$$$$ + \frac{(n-1)(2n+5)^2}{16}\Phi
- \frac{n(2n+5)\alpha}{n+2} \frac{\sum f_{,i}\rho_{,i}}{(d+f)\rho} -
h'\sum (u_{,i})^2$$$$ - \left(n+ \frac{4n\alpha}{(n-1)(n+2)}\right)h
- \frac{2n^2\alpha}{n+2}\leq 0.$$ As $\alpha = \frac{(n+2)(n-3)}{2}
+ \frac{n-1}{4}$, it is easy to check that
$$\frac{4n^2\alpha^2}{(n+2)^2(n-1)}+ \frac{2n\alpha}{n+2}
= \frac{4n^2\alpha^2}{(n+2)^2(n-1)}\left(1+
\frac{(n+2)(n-1)}{2n\alpha}\right)>\frac{4n^2\alpha^2}{(n+2)(n^2-1)}.$$
Using the Schwarz's inequality we get
$$\frac{1}{2} h\frac{\sum u_{,i}\rho_{,i}}{\rho}\leq
\frac{1}{2(n-1)}h^2\sum (u_{,i})^2+\frac{n-1}{8}\Phi,$$
$$\frac{n(2n+5)\alpha}{n+2}
\frac{\sum f_{,i}\rho_{,i}}{(d+f)\rho}\leq
\frac{4n^2\alpha^2}{(n+2)(n^2-1)}\frac{\sum (f_{,i})^2}{(d+f)^2}
+\frac{(2n+5)^2(n^2-1)}{16(n+2)}\Phi.$$ Note that
$\frac{1}{2(n-1)}h^2\geq h'$, we get from (6.5)
$$\frac{(n+2)(n-1)}{4}\Phi - \left(n+ \frac{4n\alpha}{(n-1)(n+2)}
\right)h - \frac{2n^2\alpha}{n+2}\leq 0.$$ It follows that
$$\mathcal{A}\leq 10\mathcal {B}.\leqno(6.6)$$

Secondly, we consider the following function
$$\tilde{F}= \exp \left\{-\frac{m}{C-u } + H \right\}
\frac{(h+2\alpha) \rho^{\alpha} }{(d + f)^{\frac{2n\alpha}{n+2}}}
$$defined on $S_u(p,C)$. Clearly, $\tilde{F}$ attains its supremum at
some interior point $q^*$ of $S_u(p,C)$. Thus, at $q^*$,
$$ - hu_{,i}+ \frac{h' u_{,i}}{h+2\alpha} + H_{,i}+ \alpha
\frac{\rho_{,i}}{\rho}- \frac{2n\alpha}{n+2}\frac{f_{,i}}{d+f}=
0,\leqno(6.7)$$
$$  \left(\frac{h''}{h+2\alpha}- \frac{h'^2}{(h+2\alpha)^2}
- h'\right)\sum (u_{,i})^2+\left(\frac{h'}{h+2\alpha}
-h\right)\Delta u  \leqno(6.8)$$$$+ \Delta H + \frac{n+2}{2}\alpha
\Phi - \frac{2n\alpha}{n+2}\left(\frac{\Delta f}{d+f}-\frac{\sum
(f_{,i})^2}{(d+f)^2}\right)\leq 0 $$ where
$h''=\frac{6m}{(C-u)^4}$. By (2.11) and the Schwarz inequality
$$ \sum H_{,i}^2=\sum \left(\epsilon\frac{2x_i}{(d+f)^2}-2\epsilon
\frac{\sum x_k^2 }{(d+f)^3}f_{,i}\right)^2\leq 8\epsilon H
\frac{\sum f^{ii}}{(d+f)^2}+ 8H^2\frac{\sum
(f_{,i})^2}{(d+f)^2},\leqno(6.9)$$$$ \Delta H  =
\epsilon\frac{\Delta(\sum x_k^2)}{(d+f)^2}
-4\epsilon\frac{\langle\nabla(\sum x_k^2),\nabla
f\rangle}{(d+f)^3} -2\epsilon\frac{\sum x_k^2\Delta
f}{(d+f)^3}+6\epsilon\frac{\sum x_k^2\sum
(f_{,i})^2}{(d+f)^4}\leqno(6.10)$$$$
=\frac{\epsilon}{(d+f)^2}\left[2\sum f^{ii}+\frac{n+2}{2} \langle
\nabla\log \rho,\nabla(\sum x_k^2)\rangle-4
\frac{\langle\nabla(\sum x_k^2),\nabla f\rangle}{d+f}\right]$$$$ +
6\epsilon\frac{\sum x_k^2\sum (f_{,i})^2}{(d+f)^4}-2n\epsilon
\frac{\sum x_k^2}{(d+f)^3}-(n+2)\epsilon \frac{\sum x_k^2 \langle
\nabla\log \rho,\nabla f\rangle}{(d+f)^3}$$$$\geq
\frac{\epsilon}{(d+f)^2}\sum f^{ii} -27H \frac{\sum
(f_{,i})^2}{(d+f)^2}-\frac{3(n+2)^2}{4}H\Phi-2nH. $$ Note that
$\frac{(n+2)^2}{2}>\alpha>\frac{n+2}{2}$ and
$$\Phi=\frac{1}{\alpha^2}\sum\left(-  hu_{,i}+ \frac{h' u_{,i}}{h+2\alpha}
+ H_{,i}- \frac{2n\alpha}{n+2}\frac{f_{,i}}{d+f}\right)^2$$
$$\geq\frac{1}{2\alpha^2}\sum\left(-  hu_{,i}+ \frac{h' u_{,i}}{h+2\alpha}
-\frac{2n\alpha}{n+2}\frac{f_{,i}}{d+f}\right)^2-\frac{1}{\alpha^2}\sum
(H_{,i})^2$$
$$\geq\frac{h^2}{4\alpha^2}\sum (u_{,i})^2+\frac{n^2 }{(n+2)^2 }
\frac{\sum (f_{,i})^2}{(d+f)^2} - \frac{1}{2\alpha^2}
\frac{h'^2\sum (u_{,i})^2}{(h+2\alpha)^2}-\frac{1}{\alpha^2}\sum
(H_{,i})^2-4h,$$ where we use the fact(6.4). Inserting (2.12),
(2.13), (6.9) and (6.10) into (6.8) and using the Scwartz
inequality we have $$\frac{\epsilon}{2} \frac{\sum
f^{ii}}{(d+f)^2}- a_0\Phi -a_1 h-3n\alpha\leq 0\leqno(6.11)$$ for
some constant $a_0>0$, $a_1>0$ depending only on $n$. Since $\sum
f^{ii}\geq n\rho^{\frac{n+2}{n}},$ we get
$$ \frac{\rho^{\frac{n+2}{n}}}{(d+f)^2}\leq \frac{a_0}
{\epsilon}\Phi +\frac{2a_1}{\epsilon} h+ \frac{6\alpha}{\epsilon}
.\leqno(6.12)$$ It follows that
$$ \mathcal {B}^{1+ \frac{n+2}{n\alpha}} \leq a_2\mathcal
{A}+a_3\mathcal {B},\leqno(6.13)$$ for some positive constants
$a_2$ and $a_3$, where we used the fact that
$\exp\left\{-\frac{m}{C-u}\beta\right\}h^\gamma$ has a universal
upper bounded for any $\beta >0,\ \gamma>0$. By (6.6), we have
$$\mathcal {B}\leq d_1, \quad \mathcal {A}\leq d_1\leqno(6.14)$$
for some $d_1$ depending only on $C$, $n$ and $b_1$. Thus the
proof of Lemma 6.1 is complete. \quad $\Box$ \vskip 0.1in
\noindent

In the following we estimate $\sum u_{ii}$. To this end we first
derive a general formula which we need later.

\vskip 0.1in \noindent

{\bf Lemma 6.2} {\it Let $u(\xi)$ be a smooth strictly convex
function defined in $\Omega \subset \mathbb{R}^n$. Assume that
$$\inf_{\Omega}u=0,\;\;u|_{\partial \Omega} =C.$$
Consider the function
$$F = \exp \left \{-\frac{m}{C-u} + H \right \}Q\|\nabla K\|^2,
\leqno(6.15)$$ where $Q>0$, $H>0$ and $K$ are smooth functions
defined on $\overline{\Omega}$. F attains its supremum at an
interior point $p^*$. We choose a local orthonormal frame field on
$M$ such that, at $p^*,K_{,1}=\|\nabla K\|, K_{,i}=0$, for all
$i>1$. Then at the point $p^*$ we have the following estimates
$$2\left(\frac{1}{n-1}-\delta-1\right)(K_{,11})^2 +
2\sum K_{,j}(\Delta K)_{,j}\leqno(6.16) $$$$+ 2(1-\delta)\sum
A_{ml1}^2(K_{,1})^2- \frac{(n+2)^2}{8\delta}\Phi (K_{,1})^2
 - \frac{2}{\delta(n-1)^2}(\Delta K)^2  $$$$+
\left[-h'\sum (u_{,i})^2 - h\Delta u + \Delta H + \frac{\Delta
Q}{Q} - \frac{\sum (Q_{,i})^2}{Q^2}\right] (K_{,1})^2\leq 0,$$ for
any small positive number $\delta $.} \vskip 0.1in \noindent

{\bf Proof. } We can assume that $\|\nabla K\|(p^*)> 0$.  Then, at
$p^*$, $$ F_{,i} =0,\leqno(6.17)$$$$
 \sum F_{,ii} \leq 0.\leqno(6.18) $$   By
calculating both expressions (6.17) and (6.18) explicitly, we have
$$\left(-hu_{,i} + H_{,i} + \frac{Q_{,i}}{Q}\right)\sum (K_{,j})^2 +
2\sum K_{,j}K_{,ji} =0,\leqno(6.19)$$
$$2\sum (K_{,ij})^2 + 2\sum K_{,j}K_{,jii}+ 2\sum
\left(-hu_{,i} + H_{,i} + \frac{Q_{,i}}{Q}\right) K_{,j}K_{,ji}
\leqno(6.20)$$$$ + \left[-h'\sum (u_{,i})^2 - h\Delta u + \Delta H
+ \frac{\Delta Q}{Q} - \frac{\sum (Q_{,i})^2}{Q^2}\right]
(K_{,1})^2\leq 0.$$ Let us simplify (6.20). From (6.19)
$$ 2K_{,1i} = \left(hu_{,i} - H_{,i}- \frac{Q_{,i}}{Q}\right)K_{,1}
.\leqno(6.21)$$ Applying the Schwarz inequality yields
$$ 2\sum (K_{,ij})^2  \geq 2(K_{,11})^2 +\frac{2}{n-1}(\Delta
K - K_{,11})^2 + 4\sum_{i>1}(K_{,1i})^2\leqno(6.22)$$$$  \geq
2\left(\frac{n}{n-1}-\delta \right)(K_{,11})^2 + 4\sum_{j>1}
(K_{,1j})^2 - \frac{2}{\delta(n-1)^2}(\Delta K)^2 $$ for any
$\delta> 0$. Inserting (6.21) and (6.22) into (6.20) we get
$$2\left(\frac{1}{n-1}-\delta-1\right)(K_{,11})^2 + 2\sum K_{,j}K_{,jii}
- \frac{2}{\delta(n-1)^2}(\Delta K)^2 \leqno(6.23)$$$$ +
\left[-h'\sum (u_{,i})^2 - h\Delta u + \Delta H + \frac{\Delta
Q}{Q} - \frac{\sum (Q_{,i})^2}{Q^2}\right] (K_{,1})^2\leq 0.$$  An
application of the Ricci identity shows that
$$ 2\sum K_{,j}K_{,jii}  = 2\sum K_{,j}(\Delta K)_{,j} + 2R_{11}
(K_{,1})^2\leqno(6.24) $$$$
  = 2\sum K_{,j}(\Delta K)_{,j}+ 2 \sum A_{ml1}^2(K_{,1})^2 -
(n+2) \sum A_{11k}\frac{\rho_{k}}{\rho}(K_{,1})^2$$$$\geq 2\sum
K_{,j}(\Delta K)_{,j} + 2(1-\delta)\sum A_{ml1}^2(K_{,1})^2-
\frac{(n+2)^2}{8\delta}\Phi (K_{,1})^2. $$ Consequently, inserting
(6.24) into (6.23) we get (6.16).\quad  $\Box$

\vskip 0.1in \noindent

{\bf Lemma 6.3} {\it Let $u$ be a smooth and strictly convex
function defined in $\Omega$ which satisfies the equation (1.4).
Suppose that $u$ is normalized at $p$ and the section
$\bar{S}_u(p,C)$ is compact. And assume that there are constants
$b_2\geq 0$, $d>1$ such that
$$\frac{\sum x_k^2}{(d+f)^2}\leq b_2,\quad\quad\frac{\rho^\alpha}
{(d+f)^{\frac{2n\alpha}{n+2}}}\leq b_2,\quad \quad
\frac{\rho^\alpha \Phi}{(d+f)^{\frac{2n\alpha}{n+2}}}\leq b_2$$ on
$\bar{S}_u(p,C)$. Then there is a constant $d_2>0$, depending only
on $n$, $b_2$ and $C$, such that
$$\exp \left\{-\frac{64(n-1)C}{C-u }\right\}\frac{\rho^{\alpha}\sum
u_{ii} }{(d + f)^{\frac{2n\alpha}{n+2}+2}}\leq d_2$$ on
$S_u(p,C)$, where $\alpha = \frac{(n+2)(n-3)}{2} +
\frac{n-1}{4}.$}

\vskip 0.1in \noindent {\bf Proof.} Put
$$H=\epsilon\frac{\sum x_k^2}{(d+f)^2},\quad K=x_1,\;\;\;
Q =\frac{ \rho^{\alpha}}{(d + f)^{\frac{2n\alpha}{n+2}+2}}$$ in
(6.15). Now we first calculate $2\sum K_{,j}(\Delta K)_{,j}+
2(1-\delta)\sum A_{mli}A_{mlj}K_{,i}K_{,j}.$ By (2.11) we have in
this case
$$\Delta K = \frac{n+2}{2}\langle\nabla \log \rho,
\nabla K\rangle,$$
$$2\sum K_{,j}(\Delta K)_{,j}=(n+2)\frac{\rho_{,11}}{\rho}(K_{,1})^2
- (n+2)\frac{(\rho_{,1})^2}{\rho^2}(K_{,1})^2 + (n+2)\sum
K_{,1i}K_{,1}\frac{\rho_{,i}}{\rho}$$
$$\geq  (n+2)\sum \frac{\rho_{,ij}}{\rho}K_{,i}K_{,j}-
\delta \sum (K_{,1i})^2 -\frac{(n+2)^2+1}{4\delta}\Phi
(K_{,1})^2$$ for $\delta \leq \frac{1}{4(n+2)}$. We use the
coordinates $\xi_1,..., \xi_n$ to calculate $\sum (K_{,ij})^2$ and
$\sum A_{ml1}^2(K_{,1})^2$. Note that the Levi-Civita connection
is given by $\Gamma^k_{ij} = \frac{1}{2}\sum u^{kl}u_{lij}.$ Then
$$ K_{,ij} = u_{1ij}-\frac{1}{2}\sum u_{1k}u^{kl}u_{lij}=
\frac{1}{2}u_{1ij},$$$$\sum (K_{,ij})^2 = \frac{1}{4}\sum
u^{ik}u^{jl}u_{1ij}u_{1kl},$$$$ \sum (A_{ml1})^2(K_{,1})^2  =
\frac{1}{4}\sum
u^{ik}u^{jl}u_{ijp}u_{klq}u^{pr}u_{1r}u^{qs}u_{1s}=\sum
(K_{,ij})^2.\leqno(6.25)$$ In the coordinates $x_1,...,x_n$ we
have (see (3.1))
$$\frac{\rho_{ij}}{\rho}=
\frac{\rho_{i}}{\rho}\frac{\rho_{j}}{\rho},\;\;\;\frac{\rho_{,ij}}{\rho}=
\frac{\rho_{i}}{\rho}\frac{\rho_{j}}{\rho} + \sum
A^k_{ij}\frac{\rho_k}{\rho}.$$ It follows that $$(n+2)\sum
\frac{\rho_{,ij}}{\rho}K_{,i}K_{,j}\leq \delta \sum (K_{,ij})^2+
\frac{(n+2)^2+1}{4\delta}\Phi (K_{,1})^2.\leqno(6.26)$$ Then
$$  2\sum K_{,j}(\Delta K)_{,j}+
2(1-\delta)\sum A_{mli}A_{mlj}K_{,i}K_{,j}\leqno(6.27)$$$$ \geq
(2-4\delta)\sum (K_{,ij})^2 - \frac{(n+2)^2+1}{2\delta}\Phi
(K_{,1})^2. $$  A direct calculation yields
$$\frac{\Delta Q}{Q} - \frac{\sum
(Q_{,i})^2}{Q^2} \geq -\frac{(n\alpha + n+2)(n+2)}{8}\Phi -
2n(\alpha+1).\leqno(6.28)$$ From (6.21) we obtain
$$\sum (K_{,1i})^2 =\frac{1}{4}\sum\left[hu_{,i} - \alpha
\frac{\rho_{,i}}{\rho} +\left(\frac{2n\alpha}{n+2}+2\right)
\frac{f_{,i}}{d+f}-H_{,i}\right]^2(K_{,1})^2\leqno(6.29)$$$$
\geq\frac{1}{16}\sum\left[hu_{,i}+\left(\frac{2n\alpha}{n+2}+2\right)
\frac{f_{,i}}{d+f}\right]^2(K_{,1})^2-\frac{1}{8} \alpha^2 \Phi
(K_{,1})^2-\frac{1}{4}(K_{,1})^2 \sum (H_{,i})^2 $$$$ \geq
\frac{1}{16}\left[h^2\sum (u_{,i})^2 +  \frac{4n^2\alpha^2}
{(n+2)^2}\frac{\sum
(f_{,i})^2}{(d+f)^2}\right](K_{,1})^2-\frac{1}{8} \alpha^2 \Phi
(K_{,1})^2-\frac{1}{4}(K_{,1})^2 \sum (H_{,i})^2 - a_4h (K_{,1})^2
,$$ where we used (6.4), for some positive constant $a_4$. Choose
$\delta = \frac{1}{6(n+2)}$ and $m=64(n-1)C$. Inserting (6.9),
(6.10), (6.27), (6.28) and (6.29) into (6.16) and using the
Schwarz inequality we get
$$\frac{\epsilon}{2}\frac{\sum f^{ii}}{(d+f)^2} - a_5\Phi - a_6
h - a_7\leq 0, \leqno(6.30)$$ In the above $a_4--a_7$ denote
constants depending only on $n$. Note that
$$\sum f^{ii}\geq u_{11}=(K_{,1})^2.$$ It follows that
$$\exp \left\{-\frac{m}{C-u }\right\}\frac{ \rho^{\alpha} u_{11} }{(d
+ f)^{\frac{2n\alpha}{n+2}+2}}\leq d_2$$ for some constant $d_2$
depending only on $n,\ b_2$ and $C$. Similar inequalities for
$u_{ii}$ remain true. Thus the proof of Lemma 6.3 is complete. \quad
$\Box$

\section*{\S7. Proof of Main Theorem }
\vskip 0.1in \noindent

Let $u(\xi_1,...,\xi_n)$ be a locally strongly convex function
defined on whole $\mathbb{R}^n$ such that its Legendre function $f$
satisfying
$$\frac{\partial^{2}}{\partial x_{i}\partial x_{j}}
\left(\log \det\left(f_{kl}\right)\right) = 0.\leqno(7.1)$$ Let
$p\in \mathbb{R}^n$ be any point. By a coordinate translation
transformation and by subtracting a linear function we may suppose
that $u$ satisfying
$$u(\xi)\geq u(p)=0,\ \ \ \forall \xi\in \mathbb{R}^n.$$
Choose a sequence $\{C_k\}$ of positive numbers such that
$C_k\rightarrow\infty$ as $k\rightarrow\infty$. For any $C_k$ the
level set $S_u(p,C_k)=\{u(\xi)<C_k\}$ is a bounded convex domain.
Let
$$u^{(k)}(\xi)=\frac{u(\xi)}{C_k},\quad k=1,2,\dots$$
There exists the unique minimum ellipsoid $E$ of $S_u(p,C_k)$
centered at $q_k$, the center of mass of $S_u(p,C_k)$, such that
$$n^{-\frac{3}{2}} E \subset S_u(p,C_k) \subset E.$$
Let
$$T_k: \tilde{\xi_i}=\sum a_i^j \xi_j + b_i$$
be a linear transformation such that
$$T_k(q_k) = 0,\;\;\;T_k(E)=B(0,1).$$ Then
$$B(0,n^{-\frac{3}{2}}) \subset\Omega_k:=T_k( S_u(p,C_k)) \subset
 B(0,1).$$ Thus we obtain a sequence of convex functions
$$\tilde{u}^{(k)}(\tilde{\xi}):=
u^{(k)}\left(\sum b_1^j (\tilde{\xi}_j-b_j),...,\sum b_n^j
(\tilde{\xi}_j-b_j)\right)$$ where $(b_i^j)= (a_i^j)^{-1}$.

In the following we will use the coordinates $\xi$ to denote the
$\tilde{\xi}$ and $u^{(k)}$ to denote $\tilde{u}^{(k)}$ to
simplify the notations. We may suppose by taking subsequences that
$\Omega_k$ converges to a convex domain $\Omega $ and
$u^{(k)}(\xi)$ converges to a convex function $u^\infty(\xi)$,
locally uniformly in $\Omega$. Consider the Legendre
transformation relative to $u^{(k)}$:
$$x_i=\frac{\partial u^{(k)}}{\partial \xi_i},$$ $$f^{(k)}(x_1,...,x_n)
=\sum \xi_i  \frac{\partial u^{(k)}}{\partial \xi_i}-
u^{(k)}(\xi_{1},...,\xi_{n}),\;\;\;\;(\xi_1,...,\xi_n)\in
\Omega_k.$$ Put
$\Omega^{(k)*}=\left\{(x_1,...,x_n)|x_i=\frac{\partial
u^{(k)}}{\partial \xi_i}\right\}.$ Obviously, $f^{(k)}$ satisfies
(7.1), therefore there are constants $d^{(k)}_1,...,d^{(k)}_n,
d^{(k)}_0$ such that
$$\det\left(\frac{\partial^{2}f^{(k)}}{\partial x_{i}\partial
x_{j}}\right) = \exp \left\{\sum d^{(k)}_i x_i +
d^{(k)}_0\right\}.\leqno(7.2)$$ We use Lemmas 5.1, 6.1 and 6.3 for
each $u^{(k)}$ with $C=1$ to get the following uniform estimates
$$\frac{\rho^{(k)}}{(d+f^{(k)})^{\frac{2n}{n+2}}}\leq  d_3,\;\;\;
\frac{\rho^{(k)\alpha}\Phi^{(k)}}{(d+f^{(k)})^{\frac{2n\alpha}{n+2}}}
\leq d_3,\quad \frac{\rho^{(k)\alpha}\sum
u_{ii}^{(k)}}{(d+f^{(k)})^{\frac{2n\alpha}{n+2}+2}}\leq d_3$$ on
$S_{u^{(k)}}(T^k(p),\frac{1}{2})$ for some constant $d_3>0$, where
$\alpha = \frac{(n+2)(n-3)}{2} + \frac{n-1}{4}.$\vskip 0.1in
\noindent Let $B_R(0)$ be a Euclidean ball such that
$S_{u^{(k)}}(T^k(p),\frac{1}{2})  \subset B_{R/2}(0)$, for all
$k$. The comparison theorem for the normal mapping (see[G] or
[L-J-3]) yields
$$B^*_{r}(0)\subset\Omega^{(k)*}$$
for every $k$, where $r = \frac{1}{2R}$ and $B^*_{r}(0) = \{x |
x_1^2 +...+ x_n^2 \leq r^2\}.$ Note that $u^k(T^k(p))=0$ and its
image under normal mapping is $(x_1,...,x_n)=0$. Restricting to
$B^*_{r}(0)$, we have
$$-R'\leq f^{(k)} = \sum \xi_i x_i - u^{(k)} \leq  R',$$ where
 $R'=\frac{1}{R}+1$. Therefore $f^{(k)}$ locally
uniformly converges to a convex function $f^{\infty}$ on
$B^*_{r}(0)$ and there are uniform estimates
$$\rho^{(k)} \leq  d_4,\;\;\;
(\rho^{(k)})^{\alpha}\Phi^{(k)} \leq d_4,\;\;\;
(\rho^{(k)})^{\alpha}\sum u^{(k)}_{ii} \leq d_4 \leqno(7.3)$$ on
$B^*_{r}(0)$ for some constant $d_4>0$. \vskip 0.1in \noindent

{\bf Lemma 7.1} {\it Let $f(x)$ be a smooth strictly convex
function defined in $B^*_{\delta}(0)$ satisfying $$-R'\leq f \leq
R'.$$ Then there exists a point $p^*\in B^*_{\delta}(0)$ such that
at $p^*$ $$\frac{1}{\rho} < \left(\frac{4R'}{\delta^2}\right)
^{\frac{n}{n+2}}2^\frac{n+1}{n+2} :=d_5.$$} \vskip 0.1in {\bf
Proof.}  If Lemma 7.1 does not hold, we would have
$$\frac{1}{\rho}\geq d_5\quad\quad \hbox{on}\ \ B^*_{\delta}(0).$$
It follows that
$$\det(f_{ij})\geq d_5^{n+2} \quad \hbox{on}\ \ B^*_{\delta}(0).$$
Define a function
$$F(x)=\left(\frac{d_5^{n+2}}{{2^{n+1}}}\right)^{\frac{1}{n}}
\left(\sum x_i^2-\delta^2\right) + 2R' \quad \hbox{on}\quad
B^*_{\delta}(0).$$ Then  $$\det (F_{ij})= \frac{d_5^{n+2}}{2}<\det
(f_{ij})  \quad \hbox{in}\ \ B^*_{\delta}(0),$$
$$F(x)\geq f(x)\quad\hbox{on}\quad \partial B^*_{\delta}(0).$$
By the comparison principle, we have
$$F(x)\geq f(x) \quad\hbox{on}\quad  B^*_{\delta}(0).$$
On the other hand, note that
$$F(0)=-\left(\frac{d_5^{n+2}}{2^{n+1}}\right)^{\frac{1}{n}}
\delta^2+2R'=-2R'<f(0).$$ This is a contradiction. \quad$\Box$

From Lemma 7.1 and (7.3), for any $B^*_{\delta}(0)$ we have a point
$p_k\in B^*_{\delta}(0)$ such that $\rho^{(k)}$,
$\frac{1}{\rho^{(k)}}$, $\Phi^{(k)}$ and $\sum u^{(k)}_{ii}$ are
uniformly bounded at $p_k$. Therefore there are constants $0<\lambda
\leq \Lambda < \infty$ independent of $k$ such that the following
estimates hold
$$\lambda < \hbox{the\;\;eigenvalues\;\; of\;\;}(f^{(k)}_{ij})(p_k)
< \Lambda.$$ Since $f^{(k)}$ satisfies (7.2),
$$\Phi^{(k)} = \frac{1}{(n+2)^2}\sum f^{(k)ij}d^{(k)}_id^{(k)}_j.$$
It follows that
$$\sum (d^{(k)}_i)^2 \leq d_6$$
for some constant $d_6>0$. Thus
$$\|\nabla\log\rho^{(k)}\|_E^2=\sum \left(\frac{\partial \log \rho^{(k)}}
{\partial x_i}\right)^2= \frac{1}{(n+2)^2}\sum (d^{(k)}_i)^2\leq
d_6,\leqno(7.4)$$ where $\|\cdot\|_E$ denotes the norm of a vector
with respect to the Euclidean metric. Then for any unit speed
geodesic starting from $p_k$,
$$\left|\frac{d\log\rho^{(k)}}{ds}\right|\leq\|\nabla\log\rho^{(k)}\|_E\leq
d_6.\leqno(7.5)$$ Thus for any $q$ we have
$$\rho^{(k)}(p_k)\exp\{-|q-p_k|d_6\} \leq \rho^{(k)} (q) \leq
\rho^{(k)}(p_k)\exp\{|q-p_k|d_6\}.\leqno(7.6)$$ In particular, we
choose $q$ be the point $x_i=0$ for all $i\geq 1$. It follows from
(7.3) that
$$\Phi^{(k)} (q) \leq d_7\leqno(7.7)$$
for some constant $d_7>0$ independent of $k$. On the other hand,
if $\Phi(p)\neq 0$, by a direct calculation yields
$$\Phi^{(k)}(q)= C_{k}\Phi(p)\rightarrow\infty,\;\;\;\hbox{as}\;\;\;
k\rightarrow \infty.$$ This contradicts to (7.7). Thus $$\Phi
(p)=0.$$ Since $p$ is arbitrary we conclude that $\Phi = 0$
everywhere. Consequently
$$\det\left(\frac{\partial^{2}u}{\partial \xi_{i}\partial \xi_{j}}
\right)=const.>0.$$This means that $M$ is an affine complete
parabolic affine hypersphere. By the J-C-P Theorem we conclude that
$M$ must be elliptic paraboloid. This complete the proof of the Main
Theorem. \quad $\Box$

\vskip 0.5in \noindent
An-Min Li       \hfill                              Ruiwei Xu \\
Department of Mathematics   \hfill                  Department of Mathematics \\
Sichuan University           \hfill                       Sichuan University \\
Chengdu, Sichuan              \hfill                             Chengdu, Sichuan \\
P.R.China                  \hfill                   P.R.China \\
e-mail:math-li$@$yahoo.com.cn\hfill e-mail:xuruiwei$@$yahoo.com.cn


\vskip 0.1in \noindent
\begin{thebibliography}{L3}


\bibitem[Ca]{1} E. Calabi: Improper Affine Hyperspheres of Convex
Type and a Generalization of a Theorem by K. J\"{o}rgens. Michigan
Math. J., 5(1958), 105-126.


\bibitem[C-L]{2} L. Caffarelli,  Y.Y. Li: An Extension to a Theorem of
J\"{o}rgens, Calabi, and Pogorelov. Comm. Pure Appl. Math., Vol.
LVI,(2003), 549-583.

\bibitem[C-Y]{3} S.Y. Cheng, S.T. Yau: Complete Affine Hypersurfaces. I.
The Completeness of Affine Metrics. Comm. Pure Appl. Math., 39
(1986), no. 6, 839-866.

\bibitem[C-Y-1]{4} S.Y. Cheng, S.T. Yau: On the Real Monge-Amp\`{e}re Equation
and Affine Flat Structure.  Proceedings of the 1980 Beijing
Symposium Differential Geometry and Differential Equations,
Vol.1,2,3,(Beijing), 339-370. Science Press, 1982.

\bibitem[G]{5}  C.E.  Guti\'{e}rrez: The Monge-Amp\`{e}re Equation,
Birkh\"{a}user Boston, 2001.

\bibitem[J]{6}  K. J\"{o}rgens:  \"{U}ber die L\"{o}sungen der
 Differentialgleichung $rt- s^2 = 1$. Math. Ann., 127 (1954), 130-134.

\bibitem[J-L]{7}F. Jia, A.-M. Li: Complete K\"{a}hler Affine
Manifolds. Preprint.

\bibitem[L-J-1]{10}A.-M. Li, F. Jia: Affine Bernstein Problem on Affine
Maximal Surfaces. Sichuan Daxue Xuebao, 36(1999), no.6, 1141-1143.

\bibitem[L-J-2]{10}A.-M. Li, F. Jia: The Bernstein Property of
Some Fourth Order Partial Differential Equations. Preprint.

\bibitem[L-J-3]{8}A.-M. Li, F. Jia: Euclidean  Complete Affine
Surfaces with Constant Affine Mean Curvature. Ann. Global Anal.
Geom., 23(2003), 283-304.

\bibitem [Sh]{10} H. Shima: The Geometry of Hessian
Structures, World-Scientific, 2007.

\bibitem [T-W]{11} N. Trudinger, X. Wang: The Bernstein Problem for Affine Maximal
Hypersurfaces. Invent. Math., 140(2000), 399-422.

\bibitem[P]{11} A.V. Pogorelov: On the Improper Convex Affine
Hyperspheres. Geom. Dedicata, 1 (1972), no. 1, 33-46.

\bibitem[P-1]{12} A.V. Pogorelov: The Minkowski Multidimensional
Problem. John Wiley \& Sons, 1978.

\end{thebibliography}
\end{document}